\documentclass[12pt]{article}

\usepackage{amsfonts}
\usepackage{amsthm}
\usepackage{euscript}

\usepackage[colorlinks, citecolor=black, urlcolor=blue, linkcolor=black]{hyperref}

\newtheorem{theorem}{Theorem}[section]

\newtheorem{remark}{Remark}[section]

\makeatletter
\@addtoreset{equation}{section}
\makeatother

\usepackage{fancyhdr}
\begin{document}

\thispagestyle{plain} 

\sloppy

\begin{center}
{\bf \Large\bf Spectral density in a Moszynski's class of Jacobi matrices.}
\\ \bigskip
{\large Ianovich E.A.}
\\ \bigskip
   {\it Saint Petersburg, Russia}
\\ \bigskip
   {\it\small E-mail: eduard@yanovich.spb.ru}
\end{center}

\begin{abstract}
In this paper it is considered a spectral density for a class of Jacobi matrices with absolutely continuous spectrum that was examined first by Moszynski. It is shown that the corresponding spectral density is equivalent to the positive continuous function everywhere except maybe the point $x=0$.
\end{abstract}

\section{Introduction}

Let us consider Moszynski's class of Jacobi matrices
$$
  \left(\begin{array}{ccccc} 0&b_0&0&0&\ldots\\
                      b_0&0&b_1&0&\ldots\\
                      0&b_1&0&b_2&\ldots\\
                      0&0&b_2&0&\ldots\\
                      \vdots&\vdots&\vdots&\vdots&\ddots
  \end{array}\right),
$$
where $b_0>0$, $b_{2k-1}=b_{2k}=k^\alpha\,,\:k=1,2,3,\ldots\:$; $\alpha\in(0;1)$. For given $\alpha$ this matrix describes unbounded self-adjoint operator $A$ with absolutely continuous spectrum covering the whole real axis: $\sigma_{ac}(A)=\mathbb{R}$~\cite{1}.
Standard three-term recurrence relations for this matrix have the form
\begin{equation}
\label{recurrent_relations}
b_{n-1}\,y_{n-1}+b_n\,y_{n+1}=x\,y_n\,,\quad (n=1,2,3,\ldots)
\end{equation}
and define so called generalized eigenvectors $\{y_n\}_{n=0}^\infty$. In particular, 1st type orthogonal polynomials $P_n(x)$ are solutions of~(\ref{recurrent_relations}) with initial conditions $P_0(x)=1$, $P_1(x)=x/b_0$. The relations~(\ref{recurrent_relations}) can be written down in the matrix form
$$
v_{n+1}=B_n\,v_n\,,\quad n=1,2,\ldots\,,
$$
where
$$
v_n=\left(\begin{array}{c}
                     y_{n-1} \\
                     y_n \\
                     \end{array}\right)\,,\quad
B_n=\left(
        \begin{array}{cc}
              0 & 1 \\
              -\frac{b_{n-1}}{b_n} & \frac{x}{b_n} \\
            \end{array}
          \right).
$$
In works~\cite{2,3,4,5,6} it was shown that if there exists the limit
\begin{equation}
\label{approximation}
\lim\limits_{n\to\infty}\left(b_nP^2_n(x)-b_{n-1}P_{n-1}(x)P_{n+1}(x)\right)=\lim\limits_{n\to\infty}\Delta_n(x)=\Delta(x)\,,
\end{equation}
then the corresponding spectral density of the operator $A$ equals to
$$
f(x)=\frac{1}{\pi\Delta(x)}.
$$
In particular, if the weights $b_n$ satisfy the conditions
\begin{enumerate}
\label{conditions}
\item $\displaystyle \lim b_n=+\infty$,

\item $\displaystyle \lim\frac{b_n}{b_{n+1}}=1$,

\item $\displaystyle\left\{\frac{b_{n-1}}{b_n}-\frac{b_{n-2}}{b_{n-1}}\right\}
\in l_1$,

\item $\displaystyle\left\{\frac{1}{b_n}-\frac{1}{b_{n-1}}\right\}\in l_1$\,,
\end{enumerate}
the limit function $\Delta(x)$ exists and is positive everywhere so that $f(x)\in C(-\infty,+\infty)$. In our case the third condition failures.

\section{Asymptotics of generalized eigenvectors and spectral density}

To obtain the asymptotics of orthogonal polynomials we will use the following discrete analogue of Levinson-type theorem~\cite{7}:
\begin{theorem}
\label{Levinson}
Let
$$
u_{n+1}=(A_n+V_n)\,u_n\,,\quad n=1,2,\ldots
$$
be a system of recurrence relations, where $A_n,V_n$ are $2\times2$ matrices and $u_n=(u_n^{(1)},u_n^{(2)})^{\top}$. Suppose the matrices $A_n$ are diagonalizable and not degenerate for all $n$ sufficiently large: $A_n=S_n\Lambda_n S_n^{-1}$, $\Lambda_n=diag(\lambda_1(n);\lambda_2(n))$. Let the eigenvalues of $A_n$ satisfy for all sufficiently large $n$ the condition: $|\lambda_1(n)|\le|\lambda_2(n)|$ $($or $|\lambda_1(n)|\ge|\lambda_2(n)|$$)$. Let the norms $\|S_n\|$ be bounded and
$$
S_{n+1}^{-1}S_n=E+R_n\,,\quad \sum_{n=p}^\infty\frac{\|R_n\Lambda_n\|}{|\lambda_i(n)|}<+\infty\,,\quad
\sum_{n=p}^\infty\frac{\|S_{n+1}^{-1}V_nS_n\|}{|\lambda_i(n)|}<+\infty\,,\,i=1,2\,,
$$
where $p$ is large enough, $E$ is identity matrix. Then there exist solutions of recurrence relations satisfying the following asymptotic formulas
$$
u^{(i)}_n=\left(\prod_{j=p}^{n-1}\lambda_i(j)\right)(e_n^{(i)}+\circ(1))\,,\quad n\to\infty\,,\,i=1,2\,,
$$
where $e_n^{(i)}$ are eigenvectors of $A_n$ $($$A_ne_n^{(i)}=\lambda_i(n)e_n^{(i)}$$)$.
\end{theorem}

\begin{remark}
\label{remark1}
If the matrices $A_n$ and $V_n$ depend on the parameter $x\in[a,b\,]$, then the uniform in $x$ fulfilling of all conditions of the theorem means that the term $\circ(1)$ in asymptotic formulas also uniformly tends to zero as $n\to\infty$.
\end{remark}

Let $V_k=0$, $u_k=v_{2k}$ and
$$
A_k=B_{2k+1}B_{2k}=\left(\begin{array}{cc} -\frac{b_{2k-1}}{b_{2k}} & \frac{x}{b_{2k}}\\
                                   -\frac{b_{2k-1}}{b_{2k}}\frac{x}{b_{2k+1}} & \frac{x^2}{b_{2k+1}b_{2k}}-\frac{b_{2k}}{b_{2k+1}}
                \end{array}\right)=
                \left(\begin{array}{cc} -1 & \frac{x}{k^\alpha}\\
                                   -\frac{x}{(k+1)^\alpha} & \frac{x^2}{(k+1)^\alpha k^\alpha}-\frac{k^\alpha}{(k+1)^\alpha}
                \end{array}\right).
$$
Then
$$
u_{k+1}=A_k\,u_{k}.
$$
It is easy to check that if $x\in\mathbb{R}\setminus\{0\}$, then $\mbox{discr}\,A_k<0$ for $k>N$. Hence for all $k$ sufficiently large the matrices $A_k$ have two complex conjugate eigenvalues and it is easy to check that
$$
A_k=S_k\Lambda_kS_k^{-1}\,,\quad (k>N\,,\,x\ne0)\,,
$$
where
$$
\Lambda_k=\left(
            \begin{array}{cc}
              \lambda_k^+ & 0 \\
              0 & \lambda_k^- \\
            \end{array}
          \right)\,,\quad
          \lambda_k^{\pm}=\frac{1}{2}\left(\omega_k\pm i\sqrt{4\frac{k^\alpha}{(k+1)^\alpha}-\omega_k^2}\,\right)\,,\:
\omega_k=\frac{x^2}{((k+1)k)^\alpha}-\frac{k^\alpha}{(k+1)^\alpha}-1\,,
$$
$$
S_k=\left(
            \begin{array}{cc}
              1 & 1 \\
              \frac{k^\alpha(\lambda_k^++1)}{x} & \frac{k^\alpha(\lambda_k^-+1)}{x} \\
            \end{array}
          \right)\,,
$$
so that
$$
(\lambda_k^{\pm}+1)\sim\frac{\alpha}{2k}+\frac{x^2}{2k^{2\alpha}}\pm i\frac{|x|}{k^\alpha}\,,\quad |\lambda_k^{\pm}|=\sqrt{\frac{k^\alpha}{(k+1)^\alpha}}\to 1\,,\,k\to\infty,
$$
$$
S_k\to\left(
            \begin{array}{cc}
              1 & 1 \\
              i\frac{|x|}{x} & -i\frac{|x|}{x} \\
            \end{array}
          \right)\,,\,k\to\infty.
$$
Therefore the norms $\|S_k\|$ are bounded. Next, we have
$$
S_{k+1}^{-1}S_k=\frac{1}{k^\alpha(\lambda_k^--\lambda_k^+)}\cdot
$$
$$
        \cdot\left(
            \begin{array}{cc}
              (k+1)^\alpha(\lambda_{k+1}^-+1)-k^\alpha(\lambda_{k}^++1) & (k+1)^\alpha(\lambda_{k+1}^-+1)-k^\alpha(\lambda_{k}^-+1) \\
              k^\alpha(\lambda_{k}^++1)-(k+1)^\alpha(\lambda_{k+1}^++1) &
              k^\alpha(\lambda_{k}^-+1)-(k+1)^\alpha(\lambda_{k+1}^++1) \\
            \end{array}
          \right)=E+R_k\,,
$$
where
$$
R_k\sim\frac{\frac{\alpha}{2}(\frac{1}{(k+1)^{1-\alpha}}-\frac{1}{k^{1-\alpha}})+
              \frac{x^2}{2}(\frac{1}{(k+1)^{\alpha}}-\frac{1}{k^{\alpha}})}{-2i|x|}
        \left(
            \begin{array}{cc}
              1 & 1\\
              -1 & -1
            \end{array}
          \right).
$$
It follows from this formula that $\{\|R_k\|\}\in l_1$.

Thus all conditions of the Theorem~(\ref{Levinson}) are fulfilled, and there exist solutions $u^{\pm}_k=v^{\pm}_{2k}$ such that
$$
u^{\pm}_k=v^{\pm}_{2k}=\left(\prod_{j=p}^{k-1}\lambda^{\pm}_j\right)(e_k^{\pm}+\circ(1))\,,\quad k\to\infty\,,
$$
where
$$
e_k^{\pm}=\left(
                  \begin{array}{c}
                    1 \\
                    \frac{k^\alpha(\lambda_k^{\pm}+1)}{x} \\
                  \end{array}
                \right)\to\left(
                            \begin{array}{c}
                              1 \\
                              \pm i\frac{|x|}{x} \\
                            \end{array}
                          \right)=e^{\pm}\,,\quad k\to\infty.
$$
Let $\lambda_k^{\pm}=|\lambda_k^{\pm}|\,e^{i\phi^{\pm}_k}=\frac{k^{\alpha/2}}{(k+1)^{\alpha/2}}\,e^{i\phi^{\pm}_k}$, where
$\phi^{\pm}_k\sim\pi\mp\frac{|x|}{k^\alpha}$ as $k\to\infty$. Then
$$
v^{\pm}_{2k}=\frac{p^{\alpha/2}}{k^{\alpha/2}}\,\exp\left(i\,\sum\limits_{j=p}^{k-1}\phi^{\pm}_j\right)(e^{\pm}+\circ(1))\,,\quad k\to\infty.
$$
Using Euler's summation formula, one has
$$
\sum\limits_{j=p}^{k-1}\phi^{\pm}_j=\pi k\mp\frac{|x|\,k^{1-\alpha}}{1-\alpha}\mp c+\circ(1)\,,\quad k\to\infty,
$$
where $c$ is some constant. Hence
$$
v^{\pm}_{2k}=\frac{p^{\alpha/2}}{k^{\alpha/2}}\,\exp\left[i\left(\pi k\mp\frac{|x|\,k^{1-\alpha}}{1-\alpha}\mp c\right)\right]
(e^{\pm}+\circ(1))\,,\quad k\to\infty.
$$
It follows that
$$
y_{2k-1}^{\pm}=\frac{p^{\alpha/2}}{k^{\alpha/2}}\,\exp\left[i\left(\pi k\mp\frac{|x|\,k^{1-\alpha}}{1-\alpha}\mp c\right)\right](1+\circ(1))\,,\quad k\to\infty,
$$
$$
y_{2k}^{\pm}=\pm\,\mbox{sign}(x)\,i\,\frac{p^{\alpha/2}}{k^{\alpha/2}}\,\exp\left[i\left(\pi k\mp
\frac{|x|\,k^{1-\alpha}}{1-\alpha}\mp c\right)\right]
(1+\circ(1))\,,\quad k\to\infty.
$$

It is evident that these asymptotics belong to the two linearly independent solutions of~(\ref{recurrent_relations}). Hence 1st type polynomials $P_n(x)$ are their linear combination with constant coefficients and have the following asymptotic behavior
\begin{equation}
\label{asymptotics_P_n(x)}
\begin{array}{c}
\displaystyle P_{2k-1}(x)=(-1)^k\frac{C}{k^{\alpha/2}}\,\cos\left(\frac{|x|\,k^{1-\alpha}}{1-\alpha}+\beta\right)(1+\circ(1))\,,\quad k\to\infty, \\
\displaystyle
P_{2k}(x)=(-1)^k\,\mbox{sign}(x)\,\frac{C}{k^{\alpha/2}}\,\sin\left(\frac{|x|\,k^{1-\alpha}}{1-\alpha}+\beta\right)(1+\circ(1))\,,
\quad k\to\infty,
\end{array}
\end{equation}
where $C>0$ and $\beta\in\mathbb{R}$ are some constants (depending on $x$ generally). This formulas were obtained under the condition $x\ne 0$. From remark~(\ref{remark1}) it follows that the term $\circ(1)$ tends to zero uniformly on any set of the form $[-A;-B]\cup[B;A]$ ($0<B<A$). If $x=0$ then the matrices $A_k$ become diagonal and one can easily obtain that
\begin{equation}
\label{asymptotics_P_n(0)}
\begin{array}{c}
\displaystyle P_{2k-1}(0)=0\,, \\
\displaystyle P_{2k}(0)=(-1)^k\,\frac{b_0}{k^\alpha}.
\end{array}
\end{equation}
It follows that the asymptotics at $x=0$ is changed. It reflects the fact of appearance of eigenvalue at $x=0$ when $\alpha>1/2$~\cite{1,5}.

Substituting~(\ref{asymptotics_P_n(x)}) and~(\ref{asymptotics_P_n(0)}) in~(\ref{approximation}), one obtains
\begin{equation}
\label{limit_Delta_n(x)}
\Delta(x)=\lim\limits_{n\to\infty}\left(b_nP^2_n(x)-b_{n-1}P_{n-1}(x)P_{n+1}(x)\right)=C^2\equiv C^2(x)>0\,,\:x\ne 0\,,
\end{equation}
$$
\Delta(0)=0\,,
$$
so that
$$
f(x)=\frac{1}{\pi C^2(x)}\,,\:x\ne 0.
$$

Due to uniform vanishing of $\circ(1)$ in asymptotic formulas~(\ref{asymptotics_P_n(x)}), the limit~(\ref{limit_Delta_n(x)})
is also uniform on any set of the form $[-A;-B]\cup[B;A]$ ($0<B<A$). Hence the function $\Delta(x)=C^2(x)$ is continuous everywhere except maybe the point $x=0$ ! Thus we obtain

\begin{theorem}
Spectral density of the considered class of operators is equivalent to the positive continuous function everywhere except maybe the point $x=0$.
\end{theorem}

\section{Conclusion}

Numerical calculations show that the function $\Delta(x)=C^2(x)$ is continuous and for $x=0$ so that the spectral density $f(x)$ has this point as a singularity and $f(x)\in L_p$ for some $p\in[1;\infty)$. But still it is a supposition only. It is interesting also to find the value of $p$ depending probably of $\alpha$.

\end{document}